\documentclass[12pt]{article}
\usepackage{amsfonts,amssymb}

\font\extra=msbm10 scaled\magstep1
\def\goth #1{\hbox{{\frak #1}}}
\font\extra=msbm10
\def\bbb #1{\hbox{{\extra #1}}}
\def\C{\mathbb{C}}
\def\Z{\mathbb{Z}}
\def\D{\mathbb{D}}
\def\R{\mathbb{R}}
\def\l{\lambda}
\def\g{\mathfrak g}
\def\vk{\varkappa}

\def\dsize{\displaystyle}
\def\fs{\footnotesize}
\def\nn{ \nonumber }
\def\p{ \partial }
\def\defeq {\stackrel{\mbox{\rm\small def}}{=}}

\def\bq{ \begin{equation} }
\def\eq{ \end{equation} }
\def\ben{ \begin{eqnarray} }
\def\en{ \end{eqnarray} }

\def\frac#1#2{{#1\over #2}}
\def\dfrac#1#2{{\displaystyle{#1\over#2}}}
\def\on#1#2{\mathop{\vbox{\ialign{##\crcr\noalign{\kern2pt}
$\scriptstyle{#2}$\crcr\noalign{\kern2pt\nointerlineskip}
\kern-2pt$\hfil\displaystyle{#1}\hfil$\crcr}}}\limits}

\renewcommand{\theequation}{\arabic{section}.\arabic{equation}}

\begin{document}

%%%%%%%%%%%%%%%%%%%%%%%%%%%%%%%%%%%%%%
\baselineskip=15pt
%\begin{flushright}
%Draft\\
%12/04/2003
%\end{flushright}
\vspace{1cm} \noindent {\LARGE \textbf{Compatible Lie brackets 
related to elliptic curve }} \vskip1cm \hfill
\begin{minipage}{13.5cm}
\baselineskip=15pt
{\bf
 A V Odesskii ${}^{1,} {}^{2}$ and
 V V Sokolov ${}^{1}$} \\ [2ex]
{\footnotesize
${}^1$ Landau Institute for Theoretical Physics, Moscow, Russia
\\
${}^{2}$ Manchester University, Department of Mathematics,
Manchester, UK}\\
\vskip1cm{\bf Abstract}
%\baselineskip=3D15pt

For the direct sum of several copies of $sl_n$, a family 
of Lie brackets compatible with the initial one is constructed. 
The structure constants of these brackets are expressed in terms of 
$\theta$-functions associated with an elliptic curve. 
The structure of Casimir elements for these brackets is investigated. 
A generalization of this construction to the case of vector-valued 
$\theta$-functions is presented. The brackets define a
multi-hamiltonian structure for the elliptic $sl_n$-Gaudin model. 
A different procedure for constructing compatible Lie brackets based on the 
argument shift method for quadratic Poisson brackets is discussed. 

\end{minipage}

\vglue1cm \textbf{Address}:
Landau Institute for Theoretical Physics, Kosygina 2, 119334, Moscow, Russia

\textbf{E-mail}:
odesskii@itp.ac.ru, \,  sokolov@itp.ac.ru\newpage

\centerline{\Large\bf Introduction}
\medskip
Two Lie brackets $[\cdot,\cdot]_1$ and $[\cdot,\cdot]_2$ defined on the same
finite dimensional vector space $\bf V$ are said to be compatible if 
\begin{equation} \label{pensil}
[\cdot,\cdot]_{u}=[\cdot,\cdot]_1+u [\cdot,\cdot]_2
\end{equation}
is a Lie bracket for any constant $u$.
As a matter of fact, this notion coincides with the concept of two 
compatible linear Poisson structures (see \cite{magri}). 
Indeed, the formula 
\begin{equation}\label{linpu}
\{x_i,\, x_j \}=c^k_{ij} x_k, \qquad i,j=1,\dots,N
\end{equation}
defines a Poisson bracket iff $c^k_{ij}$ are structural constants of a 
Lie algebra and the compatibility of two Poisson brackets of this form is
equivalent to the compatibility of the two corresponding Lie structures.  

The Casimir functions of the Poisson bracket $\{ \cdot,\cdot \}_{u}$ 
corresponding to (\ref{pensil}) are polynomials in $u,$  whose 
coefficients commute with respect to both Poisson brackets. This way for
constructing completely integrable Hamiltonian dynamical systems 
of compatible Poisson brackets is called the Lenard-Magri scheme 
\cite{magri, len}.  Pairs of compatible Lie brackets have
been considered in this context in \cite{bols,bolsbor}. 

However, possible applications of compatible pairs of Lie algebras 
in the integrability theory are not exhausted by this construction. 
For example, it was shown in \cite{golsok1} that the system of non-linear 
hyperbolic equations 
$$
U_x=[U,\,V]_1, \qquad  V_y=[V,\,U]_2, \qquad U,V\in {\bf V}
$$
is integrable for compatible Lie brackets. If the brackets 
$[\cdot,\cdot]_{1,2}$ coincide, this
system is just the principal chiral model.

Compatible Lie brackets also are closely related to decompositions of
infinite-dimensional Lie algebras into a vector direct sum of two subalgebras 
{\cite{reysem, skryp, golsok2}}. Furthermore, it was shown in  \cite{golsok3} 
that any  pair of compatible Lie brackets having a common quadratic Casimir 
function produces a (non-constant) solution of the classical Yang-Baxter 
equation.

The following classification problem arises: to describe all possible 
brackets $[\cdot,\cdot]_2$ on a vector space ${\bf V},$ 
compatible with a given semi-simple bracket 
$[\cdot,\cdot]_1$.
Since any semi-simple Lie algebra is rigid, the bracket (\ref{pensil}) is
isomorphic to $[\cdot,\cdot]_1$ for almost all values of the parameter 
$u$. It is well known that for the semi-simple case 
the second bracket $[\cdot,\cdot]_2$ is given by the formula 
$$
[X, \, Y]_2=[R(X), \, Y]_1+[X, \, R(Y)]_1-R([X, \, Y]_1),
$$
where $R$ is a linear operator on {\bf V}.  

Some examples 
of compatible brackets are known (see \cite{bolsbor,golsok1}). Similar to 
solutions of classical Yang-Baxter equations (\cite{beldrin}), 
these examples are (in some sense) rational, trigonometric or elliptic. 

In this paper, we construct a class of compatible semi-simple Lie brackets 
related to elliptic curve.   By analogy with other "elliptic" models in 
integrability theory, one can expect that a very wide class of compatible 
Lie brackets can be obtained by different degenerations of these basic 
"elliptic" pairs and by deformations of degenerate pairs, which are usually 
not so rigid as the elliptic models.

The paper is organized as follows. In Section 1 we present a 
construction of compatible elliptic  
$\oplus_{i=1}^m sl_n$-Lie brackets. The initial date for our construction is   
a pair of $\theta$-functions of order $m$  without
common zeros on an elliptic curve ${\cal E}$. 
To demonstrate the main idea, let us consider a slightly different situation 
of two compatible associative structures. 

Let $\bf V$ be the $k$-dimensional vector
space of all polynomials of degree $\le k-1$ in one variable,  
let $\mu_1$ and $\mu_2$ be given 
polynomials
of degree $k$ without common roots. It is clear that any polynomial $Z$, where 
$\hbox{deg} Z\le 2 k-1,$ can be uniquely represented in the form 
$
Z=\mu_1 P+ \mu_2 Q,
$
where $P, Q\in {\bf V}$. The explicit form of $P$ and $Q$ can be found with the
help of the Lagrange interpolation formula. Let us define two multiplications 
$\circ$ and $\star$ on ${\bf V}$ by the formula 
$$
X\,Y=\mu_1 (X\circ Y) +\mu_2 (X\star  Y), \qquad  X,Y\in {\bf V}.
$$
It can be checked that any linear combination of these two products is
associative. We consider an analog of this construction for Lie algebras 
replacing polynomials by $\theta$-functions on ${\cal E}$ with values 
in $sl_n$.

In Section 2 we investigate properties of the Lie algebra  
${\cal G}_u$ with the bracket (\ref{pensil}) constructed in Section 1.  
Since for generic $u$ the Lie algebra ${\cal G}_u$  is isomorphic to $\oplus_{i=1}^m sl_n,$ 
we know that the center of the universal enveloping algebra 
${\bf U}({\cal G}_u)$ is generated by $m (n-1)$ elements of 
${\bf U}({\cal G}_u)$. More precisely, we have $m$ generators of degree $p,$ 
where $p=2,3,\dots,n$. Usually the elements of the center are called 
the Casimir elements. 

Since  (\ref{pensil}) is linear in $u$, the Casimir elements can be
chosen to be polynomial in $u$. It turns out that for each $p$ there exists 
one Casimir element $K_{p,1}$ of degree $p-2$ in $u$, $m-2$ Casimir elements 
$K_{p,2},\dots, K_{p,m-1}$ of degree $p-1$ and 
one element $K_{p,m}$ of degree $p$. In particular, there exists one quadratic Casimir element, 
which does not depend on $u$. This element can be regarded as an invariant
bilinear form, common for both brackets $[\cdot,\cdot]_{1,2}.$ 

This picture is in accordance with  general results and conjectures by Gelfand and Zakharevich
\cite{gelzah}  about "good" bi-Hamiltonian structures that in our case states,
in particular, 
that  $\sum_{i,j} (2\hbox{deg}_u K_{ij}+1)$ should be equal to 
$\hbox{dim}\,{\cal G}_u=m (n^2-1).$

In Section 2 we find explicit formulas for of all these  Casimir 
elements $K_{ij}$. 

In Section 3 we generalize the results of Section 1 to the case of 
vector-valued $\theta$-functions or, which is the same, 
to the case of $l$-dimensional indecomposable holomorphic vector bundles 
over the elliptic curve. As the result, we get an 
$(l+1)$-dimensional vector space of pairwise compatible Lie brackets. 

In Section 4 we discuss a different way for constructing compatible linear 
Poisson brackets starting with a quadratic Poisson bracket of  
Sklyanin type. This construction is based on the argument shift method. 
We conjecture that the family of compatible brackets thus obtained  
coincides with the brackets constructed in Section 3.

\section{Compatible Lie brackets related to scalar $\theta$-functions}
\setcounter{equation}{0}
Let $\Gamma\subset \C$ be a lattice generated by 1 and $\tau,$ where 
$\hbox{Im}\,\tau > 0$. Let $m\in\mathbb{N}$. We denote by $\Theta_m(\tau)$ the 
vector space of holomorphic functions $\phi:\, \C \rightarrow \C$ such that 
$$
\phi(z+1)=\phi(z) , \qquad 
\phi(z+\tau)=(-1)^m \exp \left(-2 \pi i m \,z\right)\, \phi(z).
$$
Elements of this vector space are called $\theta$-functions of order $m$. 
Properties of $\theta$-functions are described, for example, in the Appendix to the 
review \cite{odes}.  In particular, the dimension of $\Theta_m(\tau)$ is equal to $m$. 
We fix a generator of $\Theta_1(\tau)$ and denote it by $\theta(z)$.  It is known that any 
element $f\in\Theta_m(\tau)$ has $m$ roots modulo $\Gamma$ and the sum of these roots is equal to zero 
modulo $\Gamma$. In particular, $\theta(z)$ has only one root modulo $\Gamma$ at $z=0$. 
If $x_1,...,x_m$ are all roots of $f\in\Theta_m(\tau)$, then
$f(z)=c\, \theta(z-x_1)\cdots\theta(z-x_m),$ 
where $c$ is a constant.

Let us fix relatively prime natural numbers $k$ and $n$ such that $1\le k < n$.
Let $\bf a$ and $\bf b$ be $n\times n$-matrices such 
that 
\begin{equation}\label{AB}
\mathbf a^n=\mathbf b^n=1, \qquad 
\mathbf b \, \mathbf a=\exp \left(\frac{2 \pi i k}{n}\right) 
\mathbf a \, \mathbf b. 
\end{equation}
Note that such a pair of matrices $(\mathbf a,\mathbf b)$ is a 
necessary ingredient of several ``elliptic'' constructions related to $sl_n$ 
(see \cite{beldrin, reysem2, odes1}). 
It is clear that the matrices 
$\mathbf a^{\alpha}\,\mathbf b^{\beta}$, where $\alpha,\beta=0,\dots,n-1,
\,\, (\alpha,\beta)\ne (0,0),\,$ form a basis in $sl_n.$ The commutator
relations between these matrices are given by
\begin{equation}\label{ab}
[\mathbf a^{\alpha_1}\,\mathbf b^{\beta_1},\,
\mathbf a^{\alpha_2}\,\mathbf b^{\beta_2}]=
\left[\exp \left(\frac{2 \pi i k \beta_1 \alpha_2}{n}\right)- 
\exp \left(\frac{2 \pi i k \beta_2 \alpha_1}{n}\right)
\right]\,\mathbf a^{\alpha_1 + \alpha_2} \mathbf b^{\beta_1+\beta_2} .
\end{equation}

We denote by ${\bf V}_m$ the  vector
space of holomorphic functions $f:\, \C \rightarrow sl_n$ satisfying the
following quasi-periodic conditions:
\begin{equation}\label{quasi}
f(z+1)=\mathbf a f(z) \mathbf a^{-1}, \qquad 
f(z+\tau)=(-1)^m \exp \left(-2 \pi i m \,z\right)\,
\mathbf b f(z) \mathbf b^{-1}.
\end{equation}
Note that if $f_1\in{\bf V}_{m_1}$ and $f_2\in{\bf V}_{m_2}$, 
then $f_1f_2\in{\bf V}_{m_1+m_2}$. It follows from
(\ref{quasi}) that if 
$$
f(z)=\sum f_{\alpha,\beta}(z)\, \mathbf a^{\alpha}\,\mathbf b^{\beta},
$$
then 
\begin{equation}
\begin{array}{c}
\label{falbet}
\displaystyle{f_{\alpha,\beta}(z+1)=\exp \left(-\frac{2 \pi i k \,\beta}{n}\right)\,
f_{\alpha,\beta}(z)}, \\[4mm]
\displaystyle{ f_{\alpha,\beta}(z+\tau)=(-1)^m \exp \left(-2 \pi i m z+
\frac{2 \pi i k \,\alpha}{n}\right)\,
f_{\alpha,\beta}(z)}.
\end{array}
\end{equation}
These identities imply that 
\begin{equation} \label{fab}
f_{\alpha,\beta}(z)=\exp \left(-\frac{2 \pi i k \,\beta}{n} \,z\right)\,
g_{\alpha,\beta}\left(z-\frac{k \alpha}{m n}-\frac{k \beta}{m n} \tau\right), 
\end{equation}
where $g_{\alpha,\beta}(z)$ belongs to $\Theta_m(\tau).$

{\bf Lemma 1.} 
{\it Let} $\mu_1,\mu_2 \in \Theta_m(\tau)$ {\it be a pair of} $\theta$-{\it functions 
that have no common zeros. Then any element} $Z \in {\bf V}_{2 m}$ 
{\it can be uniquely represented in the form}
$$
Z=\mu_1 P+\mu_2 Q, \qquad P,Q\in {\bf V}_{m}.
$$

{\bf Proof.} Consider a linear mapping $L: {\bf V}_m\oplus {\bf V}_m 
\rightarrow {\bf V}_{2 m} $ given by the formula 
$$
L(P,Q)=\mu_1 P+\mu_2 Q.
$$
We should prove that $L$ is an isomorphism. Since 
$\hbox{dim} ({\bf V}_m\oplus {\bf V}_m) = \hbox{dim} {\bf V}_{2 m}= 2 m (n^2-1),
$
it suffices to prove that $\hbox{Ker} L=0$. Substituting  
$$P=\sum P_{\alpha,\beta}(z)\, \mathbf a^{\alpha}\,\mathbf b^{\beta},\qquad
Q=\sum Q_{\alpha,\beta}(z)\, \mathbf a^{\alpha}\,\mathbf b^{\beta}$$ 
into $L(P,Q)=0,$ we find that  $\mu_1(z) P_{\alpha,\beta}(z)+
\mu_2(z) Q_{\alpha,\beta}(z)=0$ for all $(\alpha,\beta)\ne 0$. Since $\mu_1(z)$ and 
$\mu_2(z)$ have no common roots, we see that any root of $\mu_1(z)$ is a root 
of $Q_{\alpha,\beta}(z)$. We know that $\mu_1(z) \in \Theta_m(\tau)$ has exactly 
$m$ zeros modulo $\Gamma$  and the sum of all these
zeros equals $0.$ It follows from (\ref{fab}) that if $Q_{\alpha,\beta}(z)
\not\equiv 0$, then
$Q_{\alpha,\beta}(z)$ also 
has exactly $m$ zeros, but their sum is equal to $\frac{k \alpha}{n}+
\frac{k \beta}{n} \tau.$ Hence $Q_{\alpha,\beta}(z)\equiv 0,$ 
which implies that $Q\equiv P\equiv 0.$ $\blacksquare$

Using Lemma 1, for any $f_1,f_2\in {\bf V}_m$ we define  $[f_1,\,f_2]_1$ and 
$[f_1,\,f_2]_2$ by the formula
\begin{equation}\label{combrack}
[f_1,\,f_2]=\mu_1 [f_1,\,f_2]_1+\mu_2 [f_1,\,f_2]_2.
\end{equation}

{\bf Proposition 1.}  {\it The bilinear operations} $[f_1,\,f_2]_1$ {\it and} 
$[f_1,\,f_2]_2$ {\it are compatible Lie brackets on} ${\bf V}_m$.

{\bf Proof.} It is clear that the standard bracket 
$[f_1,\,f_2]=f_1 f_2-f_2 f_1$ is a Lie bracket on the vector space 
$\oplus_{p>0} {\bf V}_p.$  Lemma 1 shows that $[\cdot,\cdot]_{1,2}$ are
well-defined brackets on ${\bf V}_m.$ Substituting (\ref{combrack}) into the  
antisymmetricity condition $[f_1, f_2]+[f_2, f_1]=0$ for the standard bracket, 
we get $\mu_1 ([f_1, f_2]_1+[f_2, f_1]_1)+\mu_2 ([f_1, f_2]_2+[f_2, f_1]_2)=0$. 
It follows from Lemma 1 that $[f_1, f_2]_i+[f_2, f_1]_i=0$ for $i=1,2.$ 
Substituting (\ref{combrack}) into the Jacobi identity for the standard bracket, we 
obtain an identity of the form $\mu_1^2 P+\mu_1 \mu_2 Q+\mu_2^2 R=0$ for some 
$P,Q,R\in {\bf V}_m$. Using the same argument as in the proof of Lemma 1, 
one can prove that $P\equiv Q\equiv
R\equiv 0$. It is easy to verify that the identities $P=0$ and $R=0$ coincide 
with the Jacobi identity for the brackets $[\cdot,\cdot]_{1}$ and 
$[\cdot,\cdot]_{2}$, respectively. The identity $Q=0$ is equivalent to the 
compatibility of the brackets $[\cdot,\cdot]_{1,2}.$ $\blacksquare$

Let $x_i(u)$ be all roots of $\mu_2(z)- u\,\mu_1(z)$:
$$
\mu_2(z)- u\,\mu_1(z)=c(u) \theta(z-x_1(u))\cdots \theta(z-x_m(u)).
$$
If $x_1(u),\dots,x_m(u)$ are distinct modulo $\Gamma$, we say that $u$ is 
regular. For brevity, we use the notation $x_j$ instead of $x_j(u).$

Now we are going to prove that the linear combination 
\begin{equation}\label{bu}
[f_1,\,f_2]_u=[f_1,\,f_2]_1+u\, [f_1,\,f_2]_2
\end{equation}
of brackets (\ref{combrack}) is isomorphic to $\oplus_{i=1}^m sl_n$ for  
all regular values of $u$.

Consider the following elements $ v_{\alpha,\beta,\gamma}(u, z)
 \in {\bf V}_m$ defined by
\begin{equation}\label{vv1} 
 v_{\alpha,\beta,\gamma}(u, z)= P_{\alpha,\beta,\gamma}(u)\, 
 g_{\alpha,\beta,\gamma}(u,z),
\end{equation}
where
\begin{equation}\label{vv2} 
 P_{\alpha,\beta,\gamma}(u)=\frac{\mu_1(x_{\gamma})}
{\theta(x_{\gamma}-x_{1})\cdots \hat\gamma \cdots 
\theta(x_{\gamma}-x_{m}) \theta(-\frac{k \alpha}{n}-\frac{k \beta}{n} 
 \tau)} 
\end{equation}
and 
\begin{equation}
\begin{array}{c}
\label{vv3} 
\displaystyle{ g_{\alpha,\beta,\gamma}(u,z)= \exp \left(-\frac{2 \pi i k \,\beta}{n} \,z\right)\,
\theta(z-x_{1})\cdots \hat\gamma \cdots 
\theta(z-x_{m})\,\times} \\[4mm]
\displaystyle{\theta\left(z-x_{\gamma}-\frac{k \alpha}{n}-\frac{k \beta}{n}\tau\right)}\, \mathbf a^{\alpha}\,\mathbf b^{\beta}
\end{array}
\end{equation}
Here the symbol $\hat \gamma$ means that the factor $\theta (z-x_{\gamma})$ 
is omitted in the product,  $\alpha,\beta=0,\dots, n-1$, where $(\alpha,\beta)\ne
(0,0)$ and $\gamma=1,\dots,m.$

{\bf Theorem 1.} {\it The elements} $v_{\alpha,\beta,\gamma}(u)$ {\it satisfy the
following commutator relations }
\begin{equation}\label{comm1}
[v_{\alpha_{1},\beta_{1},\gamma_{1}},\,v_{\alpha_{2},\beta_{2},\gamma_{2}}]_{u}=0
\end{equation}
{\it for} $\gamma_1 \ne \gamma_2$ {\it and} (cf. (\ref{ab})) 
\begin{equation}
\begin{array}{c}\label{comm2}
\displaystyle{ [v_{\alpha_{1},\beta_{1},\gamma},\,v_{\alpha_{2},\beta_{2},
\gamma}]_{u}=\left[\exp \left(\frac{2 \pi i k \beta_1 \alpha_2}{n}\right)- 
\exp \left(\frac{2 \pi i k \beta_2 \alpha_1}{n}\right)
\right]\, v_{\alpha_{1}+\alpha_{2},\beta_{1}+\beta_{2},\gamma} }.
\end{array}
\end{equation}

The proof of Theorem 1 is based on the following

{\bf Lemma 2.} {\it Suppose} $x_1,\dots,x_m\in \C$ {\it are distinct modulo} 
$\Gamma$ {\it and} 
$x_1+\cdots+x_m\ne \frac{k \alpha}{n}+\frac{k \beta}{n} \tau\,\,$ {\it modulo} $\Gamma$ {\it for} 
$0 \le \alpha, \beta
< n $, $(\alpha, \beta)\ne(0,0).$  {\it Then for any} $\sigma_1,\dots, \sigma_m\in sl_n$ 
{\it there exists a unique element} $f\in {\bf V}_m$ {\it such that} 
$f(x_{\delta})=
\sigma_{\delta}$ for $\delta=1,\dots,m$.

{\bf Proof of Lemma 2.} Consider a linear mapping $M: {\bf V}_m
\rightarrow \oplus_{i=1}^m sl_n$ given by the formula
$M(f)=(f(x_1),\dots,f(x_m)).$ Since 
$\hbox{dim} {\bf V}_m = \hbox{dim} \oplus_{i=1}^m sl_n= m (n^2-1),
$
it suffices to prove that $\hbox{Ker} M=0$. Suppose $M(f)=0,$ where 
$f=\sum_{\alpha,\beta} f_{\alpha,\beta}(z) \mathbf a^{\alpha} \mathbf b^{\beta}.$
Then $f_{\alpha,\beta}(x_{\delta})=0$ for all $\alpha,\beta,\delta.$ But it
follows from (\ref{fab}) that if $f_{\alpha,\beta}\not\equiv 0$, then the 
sum of all zeros for $f_{\alpha,\beta}$  is equal to 
$\frac{k \alpha}{n}+\frac{k \beta}{n} \tau.$  $\blacksquare$

{\bf Proof of Theorem 1.} The basic idea is to verify that the identities (\ref{comm1}),
(\ref{comm2}) hold if we substitute 
$z=x_{\delta},\, \delta=1,\dots,m$. Then Lemma 2 concludes the proof.
It is easy to check that 
\begin{equation}
\begin{array}{c}\label{comst}
\displaystyle{ [g_{\alpha_{1},\beta_{1},\gamma_{1}},\,g_{\alpha_{2},\beta_{2},\gamma_2}]=
\exp \left(-\frac{2 \pi i k (\beta_1+\beta_2)}{n} \,z\right)\,\theta(z-x_{1})\cdots \hat\gamma_{1} \cdots 
\theta(z-x_{m})\,
 \times} \\[4mm]
\displaystyle{\theta(z-x_{1})\cdots \hat\gamma_{2} \cdots 
\theta(z-x_{m}) \times \theta\left(z-x_{\gamma_{1}}-\frac{k \alpha_1}{n}-\frac{k \beta_1}{n}\tau\right)
\times} \\[4mm]
\displaystyle{\theta\left(z-x_{\gamma_{2}}-\frac{k \alpha_2}{n}-\frac{k \beta_2}{n}\tau\right) \, \left[\exp \left(\frac{2 \pi i k \beta_1 \alpha_2}{n}\right)- 
\exp \left(\frac{2 \pi i k \beta_2 \alpha_1}{n}\right)
\right]\,\mathbf a^{\alpha_1 + \alpha_2} \mathbf b^{\beta_1+\beta_2}. }
\end{array}
\end{equation}
Using the formula 
\begin{equation}\label{ugu}
[\cdot,\cdot]=\mu_1 [\cdot,\cdot]_u+(\mu_2-u \mu_1) [\cdot,\cdot]_2,
\end{equation}
we find that 
$$
[g_{\alpha_{1},\beta_{1},\gamma_{1}},\,g_{\alpha_{2},\beta_{2},\gamma_2}]_{u}
(x_{\delta})=0
$$
for $\gamma_1\ne \gamma_2,$ which implies (\ref{comm1}) by Lemma 2.
If $\gamma_1=\gamma_2=\gamma,$ then it follows from (\ref{comst}), (\ref{ugu}) 
that 
$$
[g_{\alpha_{1},\beta_{1},\gamma},\,g_{\alpha_{2},\beta_{2},\gamma}]_{u}
(x_{\delta})=0
$$
for $\delta\ne \gamma$. This proves that (\ref{comm2}) holds for 
$z=x_{\delta},\, \delta\ne \gamma.$ 
A simple straightforward computation shows that (\ref{comm2}) also 
holds for $z=x_{\gamma}$.
$\blacksquare$

{\bf Remark 1.} It is possible to construct 'trigonometric' and 'rational' degenerations 
of the elliptic brackets described above. Namely, in the trigonometric case one should replace $\theta(z)$ by 
$1-\exp (2\pi iz)$, the space $\Theta_m(\tau)$ by the space of functions of the form $$a_0+a_1 \exp (2\pi iz)+\dots+
a_m \exp(2\pi imz)$$ such that $a_m=(-1)^m a_0$, and the space ${\bf V}_m$  by the space of $sl(n)$-valued  functions of the form 
$$
\sum c_{\alpha,\beta} \exp \left(2 \pi i \frac{\beta}{n} z \right) \mathbf a^{\alpha}\,\mathbf b^{\beta}, \qquad 0\le \alpha\le n, \quad 
0\le \beta \le m\,n, 
$$
where $c_{\alpha,0}=(-1)^m c_{\alpha, mn}.$ In this formula we assume that
$\mathbf a^n=\mathbf b^n=1,\,\,$  
$\mathbf b \, \mathbf a=\exp \left(\frac{2 \pi i}{n}\right) 
\mathbf a \, \mathbf b$.  

In the rational case $\theta(z)$ is replaced by $z$, the space 
$\Theta_m(\tau)$  by the space of polynomials of the form 
$\sum_{i=0}^{m} c_{\alpha} z^{\alpha}, \,\, c_{m-1}=0,$ and ${\bf V}_m$ by the
space of polynomials 
$\sum_{i=0}^{m} g_{\alpha} z^{\alpha}, \,\,g_{\alpha}\in sl_n, g_{m-1}=0$.

\newpage
\section{Structure of Casimir elements}
\setcounter{equation}{0}
\subsection{Casimir elements for $sl_n$} 
Let ${\bf e}_{\alpha},\,\, \alpha=1,\dots,n^2-1,$ be a basis in $sl_n$, and let  
${\bf e}^{\alpha}$ be the dual basis with respect to the invariant form 
$<X,Y>=\hbox{tr} (XY).$ Then the Casimir elements 
$$
C_p=\sum_{1\le \alpha_1,\dots \alpha_p\le n^2-1}  
\hbox{tr} ({\bf e}^{\alpha_{1}}\cdots {\bf e}^{\alpha_{p}})\, 
{\bf e}_{\alpha_{1}}\circ \cdots \circ {\bf e}_{\alpha_{p}}, \qquad 
p=2,\dots,n, 
$$
where $\circ$ denotes the multiplication 
in the universal enveloping algebra ${\bf U}(sl_n),$
generate the center of  ${\bf U}(sl_n).$

Let us take ${\bf t}_{\alpha,\beta}=\mathbf a^{\alpha}\,\mathbf b^{\beta}
$ for a basis in 
$sl_n,$ where $\mathbf a,\mathbf b$ are defined by (\ref{AB}), 
$0\le \alpha,\beta <n$,
and $(\alpha,\beta)\ne(0,0)$.  Then, up to a common multiplicative constant, 
the dual basis is given by 
$$
{\bf t}^{\alpha,\beta}=\exp 
\left(\frac{2 \pi i k \alpha \beta}{n}\right){\bf t}_{-\alpha,-\beta}.
$$
In this basis the Casimir elements  have the form 
$$
C_p=\sum_{\D_p} \exp 
\left(\frac{2 \pi i k}{n} \sum_{1\le j_1 \le j_2 \le p} 
\alpha_{j_{1}}\beta_{j_{2}}\right) 
{\bf t}_{\alpha_{1},\beta_{1}}\circ \cdots \circ {\bf t}_{\alpha_{p},\beta_{p}},
$$
where
$$
\begin{array}{c}
\D_p=\{\,\, (\alpha_1,\dots, \alpha_p,\,\, \beta_1,\dots \beta_p)\,\, 
\vert \quad 
0\le \alpha_i,\beta_j < n, \quad (\alpha_j,\beta_j)\ne (0,0), \\[3mm] 
\alpha_1+\cdots+\alpha_p\equiv 0 \,\,\hbox{mod}\,n,\quad
\beta_1+\cdots+\beta_p\equiv 0 \,\,\hbox{mod}\,n\, \}.
\end{array}
$$
In particular, the quadratic Casimir element is given by 
$$
C_2=\sum_{\alpha,\beta} \exp 
\left(\frac{2 \pi i k \alpha \beta}{n}\right) 
{\bf t}_{\alpha,\beta}\circ {\bf t}_{-\alpha,-\beta}.
$$

\subsection{Polynomial Casimir elements for ${\cal G}_u$} 
In the previous section, we have equipped the vector space ${\bf V}_m$ with 
the Lie bracket (\ref{bu}). For generic $u$ the corresponding Lie 
algebra ${\cal G}_u$ is isomorphic to $\oplus_{i=1}^m sl_n$ . 

It was shown that the vector space  
${\bf V}_m$ is isomorphic to 
$\oplus_{\alpha,\beta} {\bf F}_{\alpha,\beta},$
where ${\bf F}_{\alpha,\beta}$ is the vector space of holomorphic
functions satisfying (\ref{falbet}). Denote by 
$
{\bf F}_{\alpha_1,\beta_1,\dots,\alpha_p,\beta_p}
$
the vector space of holomorphic functions of
variables $z_1,\dots,z_p$ satisfying the conditions 
\begin{equation}
\begin{array}{c}
\label{fnvar}
\displaystyle{f(z_1,\dots,z_{j}+1,\dots,z_p)=
\exp \left(-\frac{2 \pi i k \,\beta_{j}}{n}\right)\,
f(z_1,\dots,z_p)}, \\[4mm]
\displaystyle{ f(z_1,\dots,z_{j}+\tau,\dots, z_p )=
(-1)^m \exp \left(-2 \pi i m z_j+
\frac{2 \pi i k \,\alpha_{j}}{n}\right)\,
f(z_1,\dots,z_p)}.
\end{array}
\end{equation}
Calculating the number of free coefficients in the Fourier decomposition of 
a function satisfying (\ref{fnvar}), we find that  
$\hbox{dim}\,{\bf F}_{\alpha_1,\beta_1,\dots,\alpha_p,\beta_p}=m^p.$ In
particular, $\hbox{dim}\,{\bf F}_{\alpha,\beta}=m$. 
Hence the vector space ${\bf F}_{\alpha_1,\beta_1,\dots,\alpha_p,\beta_p}$ 
is spanned by the products $f_{1}(z_1)\cdots f_{p}(z_p),$ 
where $f_i(z)\in {\bf F}_{\alpha_i,\beta_i}$. In other words, 
${\bf F}_{\alpha_1,\beta_1,\dots,\alpha_p,\beta_p}$ is isomorphic to 
${\bf F}_{\alpha_1,\beta_1}\otimes\cdots \otimes {\bf F}_{\alpha_p,\beta_p}.$

We use functions from 
$ {\bf F}_{\alpha_1,\beta_1,\dots,\alpha_p,\beta_p} $  to
represent elements of the universal enveloping algebra 
${\bf U}({\cal G}_u)$. Namely, to the product 
$f_{1}(z_1)\cdots f_{p}(z_p)\in 
{\bf F}_{\alpha_1,\beta_1,\dots,\alpha_p,\beta_p}$ we assign the element 
$f_1(z) {\bf t}_{\alpha_{1},\beta_{1}}\bullet\cdots \bullet 
f_p(z) {\bf t}_{\alpha_{p},\beta_{p}}\in {\bf U}({\cal G}_u),$ 
where $\bullet$ is the multiplication in ${\bf U}({\cal G}_u).$ Denote the
corresponding linear mapping from 
$\oplus_{\D_p}{\bf F}_{\alpha_1,\beta_1,\dots,\alpha_p,\beta_p} $ 
to ${\bf U}({\cal G}_u)$  by $\sigma_p.$ 

Suppose $u$ is regular, that is, all roots $x_1,\dots,x_m$ of the function 
$\mu_2(z)-u \mu_1(z)$ are distinct.
Let  $v_{\alpha,\beta,\gamma}(u,z)=s_{\alpha,\beta,\gamma}(u,z)\,
{\bf t}_{\alpha,\beta}$ be the basis of ${\bf V}_m$ given by (\ref{vv1}) -
(\ref{vv3}). 

By Theorem 1, for each $\gamma$ 
the elements of this basis satisfy the same commutator relations as 
${\bf t}_{\alpha,\beta}.$ Hence the center of ${\bf U}({\cal G}_u)$ is generated
by 
\begin{equation}\label{bedcaz}
C_{\gamma,p}=\sum_{\D_p} \exp 
\left(\frac{2 \pi i k}{n} \sum_{1\le j_1 \le j_2 \le p} 
\alpha_{j_{1}}\beta_{j_{2}}\right) 
\sigma_p \left(f_{\alpha_1,\beta_1,\dots, \alpha_p,\beta_p, \gamma} \right),
\end{equation}
where $f_{\alpha_1,\beta_1,\dots, \alpha_p,\beta_p, \gamma}(z_1,\dots,z_p)=
s_{\alpha_1,\beta_1,\gamma}(u,z_1)\cdots s_{\alpha_p,\beta_p,\gamma}(u,z_p)$.
 
Now our goal is to find linear combinations of generators (\ref{bedcaz}) that 
are polynomial in $u.$  Let $W_{\alpha_1,\beta_1,\dots, \alpha_p,\beta_p}
\subset {\bf F}_{\alpha_1,\beta_1,\dots, \alpha_p,\beta_p}$ be the vector
subspace spanned by $f_{\alpha_1,\beta_1,\dots, \alpha_p,\beta_p, \gamma},$ 
where $\gamma=1,\dots,m$. The following statement gives us an inner description
of this subspace. 

{\bf Lemma 3.} {\it The vector space} 
$W_{\alpha_1,\beta_1,\dots, \alpha_p,\beta_p}$ 
{\it consists of holomorphic functions} $f(z_1,\dots,z_p)$ {\it satisfying} 
(\ref{fnvar}) 
{\it and vanishing on the surfaces} 
$\{ z_{j_{1}}=x_{\delta_1} \, \& \, z_{j_{2}}=x_{\delta_2} \}  $ {\it for
all} $j_1\ne j_2$ {\it and} $\delta_1\ne \delta_2.$

{\bf Proof.} It is clear that the functions 
$f_{\alpha_1,\beta_1,\dots, \alpha_p,\beta_p, \gamma}(z_1,\dots,z_p)$ enjoy 
these properties. On the other hand, since $\{s_{\alpha,\beta,\gamma}(z), \, 
\gamma=1,\dots,m \}$ is a basis in ${\bf F}_{\alpha,\beta},$ the products 
$s_{\alpha_1,\beta_1,\gamma_1}(z_1)\cdots s_{\alpha_p,\beta_p,\gamma_p}(z_p),
$ where $\gamma_j=1,\dots,m,\, j=1,\dots,p$, form a basis in  
${\bf F}_{\alpha_1,\beta_1,\dots, \alpha_p,\beta_p}.$ Suppose a function
\begin{equation}\label{dek}
f(z_1,\dots, z_p)=\sum a_{\gamma_1,\dots,\gamma_p} 
\, s_{\alpha_1,\beta_1,\gamma_1}(z_1) \cdots s_{\alpha_p,\beta_p,\gamma_p}(z_p)
\end{equation}
satisfies the conditions of Lemma 3. We have to prove that 
$a_{\gamma_1,\dots,\gamma_p}=0 $ if $\gamma_{j_{1}}\ne \gamma_{j_{2}}.$
To show this, it suffices to substitute $x_{\gamma_{j_{1}}}$ and 
$x_{\gamma_{j_{2}}}$ for $z_{j_{1}}$ and 
$z_{j_{2}}$ in (\ref{dek}) and to take into account formulas 
(\ref{vv1}) - (\ref{vv3}). $\blacksquare$

Let 
$W_p \subset \oplus_{\D_p} W_{\alpha_1,\beta_1,\dots,\alpha_p,\beta_p}
$
be the vector space spanned by 
$\{ \oplus_{\D_p} f_{\alpha_1,\beta_1,\dots,\alpha_p,\beta_p,\gamma},\,\,
\gamma=1,\dots,m \}. 
$
Similarly to the proof of Lemma 3, one can prove the following 

{\bf Lemma 4.} {\it The vector space} $W_p$ {\it consists of elements of the 
form }
$\oplus_{\D_p} g_{\alpha_1,\beta_1,\dots,\alpha_p,\beta_p},$ {\it where} 
$g_{\alpha_1,\beta_1,\dots,\alpha_p,\beta_p}(z_1,\dots,z_p)$ {\it satisfies 
the conditions of Lemma 3 and, in addition,}
$$
\begin{array}{c}
\displaystyle{\exp \left(\frac{2 \pi i k}{n} x_{\gamma} (\beta_1+\cdots+\beta_p)\right)\,
g_{\alpha_1,\beta_1,\dots,\alpha_p,\beta_p}(x_{\gamma},\dots,x_{\gamma})=}
\\[4mm]
\displaystyle{\exp \left(\frac{2 \pi i k}{n} x_{\gamma} (\beta_1'+\cdots+\beta_p')\right)\,
g_{\alpha_1',\beta_1',\dots,\alpha_p',\beta_p'}(x_{\gamma},\dots,x_{\gamma})}
\end{array}
$$ 
{\it for any} $\alpha_1,\beta_1,\dots,\alpha_p,\beta_p, \,
\alpha_1',\beta_1',\dots,\alpha_p',\beta_p'$ {\it and} $\gamma=1,\dots,m.$

Using Lemmas 3,4, we construct polynomials in $u$ that span $W_p$ 
for generic $u.$ 

{\bf Theorem 2.}
{\it For arbitrary} $g(z)\in \Theta_m(\tau)$, {\it put}   
\begin{equation}\label{fkaz}
\begin{array}{c}
\displaystyle{ f_{\alpha_{1},\beta_{1},\dots,\alpha_p,\beta_p}(z_1,\dots,z_p)=
\exp\left[-\frac{2 \pi i k}{n} (\beta_1
z_1+\cdots+\beta_p z_p) \right]
\sum_{1\le t \le p}   }  g(z_t)\, \displaystyle{ \theta\left(\frac{k \alpha_t}{n}+
\frac{k \beta_t}{n} \tau\right)} \\[6mm]
\displaystyle{ 
 \times \prod_{1\le j\le p,\,\,j\ne t}\frac{\theta(z_t-z_j+\frac{k \alpha_j}{n}+
\frac{k \beta_j}{n} \tau)}{\theta(z_t-z_j)}\,
\prod_{1\le j\le p,\,\,j\ne t}(\mu_2(z_j)-u\,\mu_1(z_j)).
 }
\end{array}
\end{equation}
{\it Then} $\oplus_{\D_p} f_{\alpha_{1},\beta_{1},\dots,\alpha_p,\beta_p} $ {\it belongs to} $W_p$ {\it and, therefore, the formula}
$$
\sum_{\D_p} \exp 
\left(\frac{2 \pi i k}{n} \sum_{1\le j_1 \le j_2 \le p} 
\alpha_{j_{1}}\beta_{j_{2}}\right) 
\sigma_p \left(f_{\alpha_1,\beta_1,\dots, \alpha_p,\beta_p} \right)
$$
{\it defines a Casimir element in} ${\bf U}({\cal G}_u).$

{\bf Proof.} We must prove that 
$f_{\alpha_{1},\beta_{1},\dots,\alpha_p,\beta_p}$ satisfies the assumptions of
Lemmas 3,4. Using the quasi-periodic properties of the functions 
$\mu_1(z), \mu_2 (z)\in \Theta_m(\tau)$ and $\theta(z)\in
\Theta_1(\tau)$, one can verify condition (\ref{fnvar}) by a simple computation. 
To prove that $f_{\alpha_{1},\beta_{1},\dots,\alpha_p,\beta_p}$ is 
holomorphic, one can check that the only possible pole at $z_t=z_j$ is canceled after
summation. It is clear that if we put
$z_{j_{1}}=x_{\delta_1} \,$ and  $z_{j_{2}}=x_{\delta_2}$, where 
$\delta_1\ne \delta_2,$ then all summands in (\ref{fkaz}) vanish. Thus 
the assumptions of Lemma 3 hold. The assumption of Lemma 4 can be
checked by a straightforward computation. $\blacksquare$

{\bf Remark 2.} If $g(z)\in \Theta_m(\tau)$ does not depend on $u$, 
then the function (\ref{fkaz}) is polynomial in $u$, of degree $p-1.$

{\bf Remark 3.} Since the Casimir function given by (\ref{fkaz}) is linear in 
$g(z)$, we have constructed a linear map $T: \Theta_m(\tau)\rightarrow 
\hbox{center of} \,\, {\bf U}({\cal G}_u).$

{\bf Lemma 5.} {\it The kernel of} $T$ {\it is generated by the element}
$g(z)=\mu_2(z)-u\,\mu_1(z)$ {\it and, therefore,} $\hbox{dim\, Ker}\,T=1.$ 

{\bf Proof.} It follows from (\ref{fkaz}) that
$$
\begin{array}{c}
\displaystyle{T(\mu_2(z)-u\,\mu_1(z))=
\prod_{1\le j\le p}(\mu_2(z_j)-u\,\mu_1(z_j))\,
\exp\left[-\frac{2 \pi i k}{n} (\beta_1
z_1+\cdots+\beta_p z_p) \right]\times }\\[6mm]
\displaystyle{ \sum_{1\le t \le p}   \theta\left(\frac{k \alpha_t}{n}+
\frac{k \beta_t}{n} \tau\right)} 
\displaystyle{ 
 \prod_{1\le j\le p,\,\,j\ne t}\frac{\theta(z_t-z_j+\frac{k \alpha_j}{n}+
\frac{k \beta_j}{n} \tau)}{\theta(z_t-z_j)}.
 }
\end{array}
$$
Consider the function 
$$
\frac{T(\mu_2(z)-u\,\mu_1(z))}{\prod_{1\le j\le p}(\mu_2(z_j)-u\,\mu_1(z_j))}.
$$
It can be checked that this function is holomorphic and satisfies 
(\ref{fnvar}) with $m=0$. Analyzing its Fourier decomposition, we see that 
such a function is identically zero. Hence $T(\mu_2(z)-u\,\mu_1(z))=0.$
Suppose now that $T(g(z))=0.$ Substituting a root $x$ of 
$\mu_2(z)-u\,\mu_1(z)$ for $z$ in (\ref{fkaz}), we see that $g(x)=0$. Since 
 $g(z)$ has exactly $m$ zeros $\hbox{mod}\, \Gamma,$ the function $g(z)$ is
 proportional to  $\mu_2(z)-u\,\mu_1(z).$ $\blacksquare$
 
 It follows from Lemma 5 that $T(\mu_2)$ is a polynomial of degree $p-2$ in 
 $u$. Therefore, formula (\ref{fkaz}) yields an $(m-1)$-dimensional subspace in 
 the $m$-dimensional vector space $W_p$ such that one generator of 
 this subspace is the polynomial $T(\mu_2)$ 
 of degree $p-2$ and $m-2$ generators 
 are polynomials of degree $p-1$ in $u$.
 
 In the following statement we construct a remaining generator of $W_p.$
 
 {\bf Theorem 3.} {\it Let }
$$
\begin{array}{c}
\displaystyle{ h_{\alpha_{1},\beta_{1},\dots,\alpha_p,\beta_p}(z_1,\dots,z_p)=
\sum_{1\le t \le p}   }  \displaystyle{A_t\,\times  
(\mu_2'(z_t)-u \mu_2'(z_t))\,\,
\prod_{1\le j\le p,\,\,j\ne t} (\mu_2(z_j)-u\,\mu_1(z_j)) } \,\,- \\[6mm]
\displaystyle{B\,\times \prod_{1\le j\le p} (\mu_2(z_j)-u\,\mu_1(z_j)),
 }
\end{array}
$$
{\it where} $A_t$ {\it is given by}
$$
A_t=\exp\left[-\frac{2 \pi i k}{n} (\beta_1
z_1+\cdots+\beta_p z_p)\right] \theta\left(\frac{k \alpha_t}{n}+
\frac{k \beta_t}{n} \tau\right)\displaystyle{ 
 \prod_{1\le j\le p,\,\,j\ne t}\frac{\theta(z_t-z_j+\frac{k \alpha_j}{n}+
\frac{k \beta_j}{n} \tau)}{\theta(z_t-z_j)}\,}
$$
{\it and} $B$ {\it is defined by the formula}
$$
B=\frac{m}{n} \exp\left[-\frac{2 \pi i k}{n} (\beta_1
z_1+\cdots+\beta_p z_p)\right] \left(B_1+B_2\right),
$$
{\it where}
$$
\begin{array}{c}
\displaystyle{B_1=\
\sum_{1\le t \le p,\,
1\le j \le p,\, j\ne t}  
   \theta \left(\frac{k \alpha_t}{n}+
\frac{k \beta_t}{n} \tau\right) \frac{\theta' (z_t-z_j+\frac{k \alpha_j}{n}+
\frac{k \beta_j}{n} \tau)}{\theta(z_t-z_j)}\,\times}  \\[6mm] 
\displaystyle{ \prod_{1\le l\le p,\,\,l\ne j,t} 
\frac{\theta(z_t-z_l+\frac{k \alpha_l}{n}+
\frac{k \beta_l}{n} \tau)}{\theta(z_t-z_l)}
 }
\end{array}
$$

$$
B_2=\sum_{1\le t \le p}    \theta' \left(\frac{k \alpha_t}{n}+
\frac{k \beta_t}{n} \tau\right)
 \prod_{1\le l\le p,\,\,l\ne t}\frac{\theta(z_t-z_l+\frac{k \alpha_l}{n}+
\frac{k \beta_l}{n} \tau)}{\theta(z_t-z_l)}.
$$
{\it Then the formula}
$$
\sum_{\D_p} \exp 
\left(\frac{2 \pi i k}{n} \sum_{1\le j_1 \le j_2 \le p} 
\alpha_{j_{1}}\beta_{j_{2}}\right) 
\sigma_p \left(h_{\alpha_1,\beta_1,\dots, \alpha_p,\beta_p} \right)
$$
{\it defines a Casimir element in} ${\bf U}({\cal G}_u).$

{\bf Proof} is similar to the proof of Theorem 2.

{\bf Remark 4.} It is clear that the Casimir element constructed in Theorem 3  
is polynomial in $u$, of degree $p.$

\section{Families of compatible Lie brackets associated with vector
$\theta$-functions }
\setcounter{equation}{0}

In this section we generalize the construction in Section 1 replacing the usual $\theta$-functions 
by vector-valued $\theta$-functions. All proofs are similar to those in 
Section 1. 

Let $\Gamma\subset \C$ be a lattice spanned by 1 and $\tau,$ where 
$\hbox{Im}\,\tau > 0$. Our general construction will depend on $d, l ,m \in \mathbb{N}$ such that $1\le l < m$
and $m,l$ are relatively prime. Denote by $ V\Theta^d_{m/l}$ the vector space consisting of holomorphic 
functions $f: \, \C^{l+1}\rightarrow \C$ of variables $z$, $x_0,\dots,x_{l-1},$ possessing the following properties:
\begin{itemize}
\item $f(z, x_0,\dots, x_{l-1})$ is a homogeneous polynomial of degree $d$ in variables $x_0,\dots, x_{l-1}$, 
\item   $$
f(z+1, x_0,\dots,x_{l-1})=f(z, 
{\bf p} (x_0),\dots,{\bf p} (x_{l-1})),
$$
\item  $$ \displaystyle f(z+\tau, x_0,\dots,x_{l-1})=\exp\left[-2 \pi i \left(\frac{m}{l} z+
\frac{m-l-1}{2 l}  \right)\,d \right]\, f(z, 
{\bf q} (x_0),\dots,{\bf q} (x_{l-1})),
$$
\end{itemize}
where 
$$
{\bf p} (x_\alpha)= \exp\left(-2 \pi i \frac{m}{l} \alpha\right)\, x_{\alpha}, 
\qquad {\bf q} (x_\alpha)=  x_{\alpha+1}, \qquad \alpha\in \Z/l\Z.
$$

{\bf Lemma 6.}  $$ \hbox{dim} \, V\Theta^d_{m/l}=m\,
\frac{(l+1)\cdots (l+d-1)}{(d-1)!}. $$

It follows from this formula that in the case $d=1,$  which is of most
importance for us, $\hbox{dim} V\Theta^1_{m/l}=m.$

{\bf Remark 5.}  Our space $ V\Theta^d_{m/l}$ is a space of holomorphic sections 
of an indecomposable vector bundle of degree $m$ and rank $l$ on the elliptic curve.  
The classification of holomorphic vector bundles on elliptic curves was obtained 
in the paper \cite{atia}.

Let $1\le k < n$ and let $k,n$ be relatively prime. Denote by ${\bf V}^d_{m,l}$ 
the vector space of all holomorphic functions $f: \C^{l+1} \rightarrow sl_n$ 
such that  
\begin{itemize}
\item $f(z, x_0,\dots, x_{l-1})$ is a homogeneous polynomial of degree $d$ in variables $x_0,\dots, x_{l-1}$, 
\item   $$
f(z+1, x_0,\dots,x_{l-1})={\bf a}\,f(z, 
{\bf p} (x_0),\dots,{\bf p} (x_{l-1})) \,{\bf a^{-1}}
$$
\item  $$  f(z+\tau, x_0,\dots,x_{l-1})=\exp\left[-2 \pi i \left(\frac{m}{l} z+
\frac{m-l-1}{2 l}  \right)\,d \right]\, {\bf b}\, f(z, 
{\bf q} (x_0),\dots,{\bf q} (x_{l-1}))\,{\bf b}^{-1},
$$
\end{itemize}
where ${\bf a}$ and ${\bf b}$ satisfy (\ref{AB}).

{\bf Lemma 7.} {\it Suppose} $\mu_1,\dots, 
\mu_{l+1}\in V\Theta^1_{m/l}$  {\it have no common
zeros for} $(x_0,\dots,x_{l-1})\ne 0.$ {\it Then any element} 
$Z\in {\bf V}^2_{m,l}$ 
{\it can be uniquely represented in the form} 
$$
Z=\mu_1 P_1+\cdots+\mu_{l+1} P_{l+1}, \qquad  P_i\in {\bf V}^1_{m,l}.
$$

It is clear that if $f,g\in {\bf V}^1_{m,l},$ then $f g-g f$ belongs to 
${\bf V}^2_{m,l}.$ Using Lemma 7, we define $[\cdot,\cdot]_1,\dots,
[\cdot,\cdot]_{l+1}$ by the formula 
$$
f_1 f_2-f_2 f_1=\mu_1 [f_1,f_2]_{1}+\cdots+\mu_{l+1} [f_1,f_2]_{l+1}, \qquad 
f_1,f_2\in {\bf V}^1_{m,l}.
$$

{\bf Proposition 2.}  {\it The bilinear operations} $[\cdot,\cdot]_1,\dots,
[\cdot,\cdot]_{l+1}$ {\it are Lie brackets on the $m$-dimen\-sional vector space} 
${\bf V}^1_{m,l}$. {\it All these Lie brackets are pairwise compatible}.

It is clear that any linear combination of brackets from Proposition 2 
is a Lie bracket. We call a $d$-dimensional vector space 
of pairwise compatible Lie brackets a $d$-{\it Lie structure}.

{\bf Remark 6.} Suppose that one of the sections $\mu_1,\dots, 
\mu_{l+1}\in V\Theta^1_{m/l}$, say, $\mu_{l+1}$, is nonzero for each $z$. In this case the subbundle 
generated by $\mu_{l+1}$ is trivial. Consider the quotient bundle modulo  
this subbundle. It has degree $m$ and 
rank $l-1$. It is clear that the $l$-Lie structure obtained from this quotient bundle is a substructure of 
our $(l+1)$-Lie structure. Counting of parameters shows that any generic $l$-Lie structure is obtained in 
this way. Therefore, any $(l+1)$-Lie structure constructed in this section is embedded into an $m$-Lie structure 
corresponding to $l=m-1$.

\section{Argument shift method for quadratic Poisson brackets}
\setcounter{equation}{0}
 
The standard argument shift method allows one to get a family of constant Poisson 
brackets compatible with any linear Poisson bracket (\ref{linpu}). Namely,  
if we perform a shift of coordinates 
$x_i \mapsto x_i+u a_i,$ where $a_i$ are 
arbitrary constants, we will have as the result an inhomogeneous linear 
bracket of the form $\{\cdot,\cdot\}_{u}=\{\cdot,\cdot\}+u 
\{\cdot,\cdot\}_1$, where the operation $\{\cdot,\cdot\}_1$ is a constant
Poisson bracket depending on the shift vector ${\bf a}=(a_1,\dots,a_N).$ 
Moreover, since the shift vector ${\bf a}$ is arbitrary, we have got an 
$N$-dimensional vector space of constant Poisson brackets such that each of these 
brackets is compatible with (\ref{linpu}) and any two of them 
are pairwise compatible.

Consider now the case of a finite-dimensional quadratic Poisson bracket.
Suppose we have a Poisson bracket of the form 
\begin{equation}\label{qua}
\{x_i,\,x_j\}=\Gamma_{i,j}^{p,q} \,x_p x_q, \qquad i,j=1,\dots,N
\end{equation}
The shift $x_i\rightarrow x_i+u a_i$ yields a Poisson bracket of the
form $\{\cdot,\cdot\}_{u}=\{\cdot,\cdot\}+u 
\{\cdot,\cdot\}_1+u^2 \{\cdot,\cdot\}_2.$
If the coefficient of $u^2$ is equal to zero, then 
this formula defines a compatible pair consisting of the quadratic bracket
(\ref{qua}) and a  
linear Poisson bracket. 
This means that the shift vector
${\bf a}$ is not arbitrary one but satisfies the following overdetermined system of algebraic
equations: 
\begin{equation}\label{admiss}
\Gamma_{i,j}^{p,q} \,a_p a_q=0,  \qquad i,j=1,\dots,N.
\end{equation}
Such a vector is said to be {\it admissible}. It is clear that the set of admissible vectors coincides 
with the set of zero-dimensional symplectic leaves of the Poisson structure (\ref{qua}). Note that 
if the set of admissible vectors contains a $p$-dimensional vector space, then shifting by vectors 
from this space we obtain $p$ compatible linear brackets and each of them is compatible with the quadratic 
bracket (\ref{qua}).

Let us apply this construction to quadratic elliptic Poisson structures (see \cite{odes}).
For most of these brackets, the system of equations (\ref{admiss}) 
has no non-trivial solutions. Nevertheless, for some important brackets 
of Sklyanin type  non-trivial admissible vectors exist. 

{\bf Example.} Consider the following  
quadratic Poisson brackets between variables $x_0,\dots,x_7$ 
(subscripts are taken modulo 8):
\begin{equation} \label{q38}
\begin{array}{l}
\{ x_i, \, x_{i+1}\}= p_{1} x_{i} x_{i+1}+k_{1} x_{i+2} x_{i+7}-
2 k_{2} x_{i+3} x_{i+6}+p_{2} x_{i+4} x_{i+5}, 
\\[3mm]
\{ x_i, \, x_{i+2}\}= p_{3} (x_{i+1}^2 - x_{i+5}^2), 
\\[3mm]
\{ x_i, \, x_{i+3}\}= p_{1} x_{i} x_{i+3}+
k_{1} x_{i+5} x_{i+6}-2 k_{2} x_{i+1} x_{i+2}+p_{2} x_{i+4} x_{i+7}, 
\\[3mm]
\{ x_i, \, x_{i+4}\}= p_{4} (x_{i+1} x_{i+3} - x_{i+5} x_{i+7}), 
\end{array}
\end{equation}
where 
$$
p_1=-\frac{1}{2} k_1^{1/2} k_2^{-1/2}(4 k_2^2+k_1^2)^{1/2}, \qquad 
p_2= k_2^{1/2} k_1^{-1/2}(4 k_2^2+k_1^2)^{1/2},
$$
$$
p_3=k_2^{1/4} k_1^{1/4}(4 k_2^2+k_1^2)^{1/4}, \qquad 
p_4=k_2^{-1/4} k_1^{-1/4}(4 k_2^2+k_1^2)^{3/4}
$$
$k_1,k_2$ are arbitrary parameters. These brackets depend on the only essential
parameter $k_1/k_2.$

Brackets (\ref{q38}) possess the following four Casimir functions 
$$
C_i=k_2 (x_i^2+x_{i+4}^2)+p_3 (x_{i+3} x_{i+5}+x_{i+1} x_{i+7})+k_1 x_{i+2}
x_{i+6}, \qquad i=0,1,2,3.
$$
The admissible vectors  are given by 
$$
{\bf a}_{\pm}=
(t_1,0,t_2, 0, \pm t_1, 0, \pm t_2, 0), \qquad {\bf b}_{\pm}=
(0,t_1,0,t_2, 0, \pm t_1, 0, \pm t_2),
$$
where $t_1,t_2$ are arbitrary parameters. We see that the admissible vectors 
form four 2-dimensional vector spaces such that $\R ^8$ is their direct sum. 

Consider the shift of coordinates defined by ${\bf a}_{+}.$ 
As the result, we get a linear bracket $\{\cdot,\cdot\}_a=t_1 \{\cdot,\cdot\}_1+
t_2 \{\cdot,\cdot\}_2.$  Hence, we obtain a pair of compatible linear 
Poisson brackets $\{\cdot,\cdot\}_{1,2}.$  For generic $t_1,t_2,$ the bracket 
$\{\cdot,\cdot\}_a$ is isomorphic to $gl_2\oplus gl_2.$ It is easy to verify
that the bracket $\{\cdot,\cdot\}_a$ has two linear Casimir functions 
$K_1=x_0+x_4$ and $K_2=x_2+x_6.$ After reducting the linear brackets to the
surface $K_1=K_2=0$, we get a pair of compatible $sl_2\oplus sl_2$ brackets. 
It is important to mention that the initial quadratic
bracket (\ref{q38}) cannot be restricted to the surface $K_1=K_2=0$ since 
$K_i$ are not Casimir functions for (\ref{q38}). One can check that the
Lenard-Magri scheme applied to the reduced brackets 
$\{\cdot,\cdot\}_{1,2}$ produces the so(4) Schottky-Manakov top.

In the paper \cite{odes}, the Poisson algebra (\ref{q38}) is denoted by $q_{8,3}(\tau)$. 
It turns out that the situation is the same for a wide class of quadratic elliptic
Poisson algebras.  

{\bf Theorem 4.} {\it For the quadratic Poisson algebras} 
$q_{mn^2,kmn-1}(\tau)$ ({\it in 
the notations of} \cite{odes}), {\it 
the set of admissible vectors is a union of} $n^2$ {\it components which are} 
$m$-{\it dimensional 
vector spaces. The space of generators of the algebra is the direct sum of 
these spaces.}

Theorem 4 can be proved  using the so-called functional realization of 
these Poisson algebras (see \cite{odfeig,odes}).
The proof will be given in another publication. The case $m=1$ was 
considered in details in \cite{olsh}.

It is clear that for any $m$-dimensional vector space of admissible vectors of the Poisson algebra 
$q_{mn^2,kmn-1}(\tau)$ one obtains after shifting by these vectors, $m$ compatible linear Poisson structures.

{\bf Conjecture 1.}  Each of the corresponding Lie algebras is isomorphic to $\oplus_{i=1}^m gl_n$. Moreover, all these Lie algebras 
have a common center. After factorization with respect to the center one obtains $m$ compatible $\oplus_{i=1}^m sl_n$ brackets. 
These $m$-Lie structures are isomorphic to the one constructed in Section 3, where $l=m-1$.

{\bf Conjecture 2.}  Each of the $(l+1)$-Lie structures constructed in Section 3 is a substructure of this $m$-Lie 
structure.

{\bf Remark 7.} The Lenard-Magri scheme applied to the pair of compatible 
linear brackets described in Section 1,2 gives rise to an integrable model 
with the $\oplus_{i=1}^m sl_n$ Poisson brackets.
Probably this integrable system is nothing but the elliptic Gaudin model 
\cite{gaud, skltak}. However, the family of integrals  
for the $\oplus_{i=1}^m sl_2$-Gaudin model considered in \cite{skltak} contains 
one parameter related to the elliptic curve plus $m-1$ additional constant
parameters. In our construction, we have $2 m - 2$ additional parameters. 
But if Conjectures 1 and 2 are true, then all these additional parameters are
inessential in the following sense. The complete set of integrals is given by 
the Casimir functions of the quadratic brackets. These integrals depend on the
elliptic curve only. Furthermore, there exist linear brackets 
$\{\cdot,\cdot\}_1,\dots, \{\cdot,\cdot\}_{m}$ that depend on the elliptic
curve only such that any linear combination of these brackets is a 
Poisson bracket as well. The integrals commute with respect to the 
whole family of these brackets. If we choose two generic brackets 
$$
\sum_{i=1}^{m} c_i \{\cdot,\cdot\}_i, \qquad \hbox{and} \qquad
\sum_{i=1}^{m} \bar c_i \{\cdot,\cdot\}_i
$$
and bring the first one to the canonical form $\oplus_{i=1}^m sl_n$ by a 
linear transformation, then the coefficients $c_i$ appear as   
parameters in the integrals and the second bracket   
becomes dependent on parameters $c_i, \bar c_i.$ 

{\bf Remark 8.} One more construction of families of compatible linear Poisson
brackets is known \cite{reysem}. It would be interesting to understand whether 
these families coincide with those described in our paper or not.

{\bf Acknowledgments.}  The paper has been written during the visit of both
authors to the Max Planck Institute (Bonn). The authors are grateful 
to MPI for hospitality and financial support. 
The authors are grateful to G. Falqui, S. Lando and M A Semenov-{T}ian-{S}hansky
for useful discussions. The research was partially supported by short visit
grant ESF-643, RFBR grant 05-01-00189, NSh 1716.2003.1 and 2044.2003.2.

\end{document}